\documentclass[12pt]{article}
\usepackage{amssymb}
\usepackage{amsmath}
\usepackage{t1enc}

\usepackage[cp1250]{inputenc}
\newtheorem{Theorem}{Theorem}
\newtheorem{Proposition}{Proposition}
\newtheorem{Lemma}{Lemma}
\newtheorem{Remark}{Remark}
\newtheorem{Corollary}{Corollary}
\begin{document}

\title{Borsuk-Ulam theorem for the loopspace\\ of a sphere}
\date{}
\author{Dariusz Miklaszewski}
\maketitle

\begin{abstract}
We study some Borsuk-Ulam type results for the loopspace of an
euclidean sphere without loops equal to their inverses . \end{abstract}

\section{Introduction.} Let $s_0 \in \mathbb{S}^n$ and
\[\Omega \mathbb{S}^n := \{\alpha : \mathbb{I} \rightarrow \mathbb{S}^n : \alpha(0) = \alpha(1) = s_0\}.\]
The group $\mathbb{Z}_2$ acts on $\Omega \mathbb{S}^n$ by
\[\alpha \mapsto \alpha^\star,~~\alpha^\star(t) := \alpha(1-t).\]
The set of fixed points of this action is denoted by
\[\Omega_{\mathrm{tf}} := \{\alpha \in \Omega \mathbb{S}^n : \alpha = \alpha^\star\},\]
where tf is an abbreviation for {\it "to and fro"}.\\ For any
mapping $f : \Omega \mathbb{S}^n \setminus \Omega_{\mathrm{tf}}
\rightarrow \mathbb{R}^k$ denote by $A_f$ the set
\[A_f := \{\alpha \in \Omega \mathbb{S}^n \setminus \Omega_{\mathrm{tf}} : f(\alpha) = f(\alpha^\star)\}.\]
Main result is the following theorem:
\begin{Theorem}
(A)~If $k<n$ then for every $f : \Omega \mathbb{S}^n \setminus
\Omega_{\mathrm{tf}} \rightarrow \mathbb{R}^k$~ $\mathrm{dim} A_f
= \infty$.\vspace{2mm}

(B)~If $k=n$ then for every $f : \Omega \mathbb{S}^n \setminus
\Omega_{\mathrm{tf}} \rightarrow \mathbb{R}^k$~ $A_f \neq
\emptyset$.
\end{Theorem}
Proof is based on the classical Borsuk-Ulam theorem \cite{B} and
on the Jaworowski-Nakaoka theorem \cite{J}, \cite{N}. The next
proposition needs the following
\begin{Lemma}
The inclusion~ $\Omega X \setminus \Omega_{\mathrm{tf}}
\hookrightarrow \Omega X$ is a homotopy equivalence for every
space $X$ with a nontrivial path from its base point.
\end{Lemma}
Denote by $w_1$ the first Stiefel-Whitney class of the double
covering $p: E \rightarrow B$ with
\[E:= \Omega \mathbb{S}^n \setminus \Omega_{\mathrm{tf}}~~\mathrm{and}~~B:= E/\mathbb{Z}_2.\]
It appears that $(w_1)^{n-1} \neq 0$ and
\begin{Proposition}
If $n \geq 3$ then $(w_1)^n =0.$
\end{Proposition}

\begin{Corollary}
There is no an equivariant mapping~ $\mathbb{S}^n \rightarrow
\Omega \mathbb{S}^n \setminus \Omega_{\mathrm{tf}}$, ($n\geq3$).
\end{Corollary}

\begin{Corollary} The proof of $A_f \neq \emptyset$ based on
the commutative diagram
\[\begin{array}{ccc}
E& \stackrel{\varphi}{\longrightarrow}& \mathbb{S}^{k-1}\\
p\downarrow& & \downarrow \gamma\\
B& \stackrel{\bar{\varphi}}{\longrightarrow}& \mathbb{RP}^{k-1}
\end{array}
\]
with $\varphi(\alpha)=
||f(\alpha)-f(\alpha^\star)||^{-1}(f(\alpha)-f(\alpha^\star))$ and
on the contradiction:
\[0= \bar{\varphi}~^\star((w_1(\gamma)^k) = (w_1)^k \neq 0\]
works for $k<n$ only, giving a weaker result than Theorem 1.
\end{Corollary}

\section{Proof of Lemma 1.}
Without loss of generality consider $X := \mathbb{S}^n$. We start
with the description of the way, how to push a loop off the set
$\Omega_{\mathrm{tf}}$. Fix two arcs $\mu$, $\nu$ of a great
circle of $\mathbb{S}^n$ going in opposite directions from $s_0$,
as their parameter inreases from $0$ to $1$. Assume that $\mu$,
$\nu$ are of the same length. For this proof it suffices that $\mu
\neq \nu$. We made stonger assumptions on $\mu$, $\nu$, which will
be essential for further considerations. Put $\delta :=
\frac{1}{4}$. For every $\alpha \in \Omega \mathbb{S}^n $ define
$\alpha_\mu \in \Omega \mathbb{S}^n \setminus
\Omega_{\mathrm{tf}}$ as the concatenation of paths:
\[\mu~\mathrm{on}~[0,\frac{\delta}{2}],~ \mu^\star ~\mathrm{on}~[\frac{\delta}{2}, \delta],~\alpha ~\mathrm{on}~[\delta,1-\delta],~
\nu~\mathrm{on}~[1-\delta,1-\frac{\delta}{2}],~
\nu^\star~\mathrm{on}~[1-\frac{\delta}{2},1].\] The same schema
defines $\alpha_{s\mu}$ for $s \in [0,1]$ with
\[\delta:= \frac{s}{4},~\mu(t):= \mu(st),~\nu(t):= \nu(st).\]
In other words,
\[\alpha_{s\mu}(t) = \left\{\begin{array}{lll}
\mu(8t)& \mathrm{for}& t \in [0,\frac{s}{8}]\\
\\
\mu(2s-8t)&\mathrm{for}& t \in [\frac{s}{8},\frac{s}{4}]\\
\\
\alpha\left(\frac{2}{2-s}(t-\frac{s}{4})\right)&\mathrm{for}& t\in
[\frac{s}{4},1-\frac{s}{4}]\\
\\
\nu(8[t+\frac{s}{4}-1])&\mathrm{for}& t \in
[1-\frac{s}{4},1-\frac{s}{8}]\\
\\
\nu(8(1-t))&\mathrm{for}& t \in [1-\frac{s}{8},1].
\end{array}\right.
\]
In particular, $\alpha_{0\mu} = \alpha$, $\alpha_{1\mu}=
\alpha_\mu$. In this way, the mapping
\[\Omega \mathbb{S}^n \ni \alpha \mapsto \alpha_\mu \in \Omega \mathbb{S}^n \setminus \Omega_{\mathrm{tf}}\]
is a homotopy inverse to the inclusion. $\Box$

\section{Embeddings of $\mathbb{S}^{n-1}$ in $\Omega
\mathbb{S}^n$.} Imagine the sphere $\mathbb{S}^n$ with its equator
$\mathbb{S}^{n-1}$, the south pole $s_0, (x_{n+1} =-1)$ and
another base point $s_1$ on the equator.\vspace{2mm}

(3.1) [Embedding $\alpha$.] If $x \in \mathbb{S}^{n-1}$ then the
loop  $\alpha(x)$ is a (great) circle in the intersection of
$\mathbb{S}^n$ with a plane parallel to the axis $OX_{n+1}$ and
going through points $s_0$, $x$. Its orientation is chosen in this
way, that $x$ is passed before $-x$.\vspace{2mm}

(3.2) [Embedding $\bar{R}\circ\beta$.] If $x \in \mathbb{S}^{n-1}$
then the loop $\beta(x)$ is a circle in the intersection of
$\mathbb{S}^n$ with a plane parallel to the axis $OX_{n+1}$ and
going through points $s_1$, $x$. Its orientation is chosen in this
way, that the south half-sphere is passed before the north
one.\vspace{2mm}

Moving the base point from $s_0$ to $s_1$ (along the path
$s_\lambda$ on the meridian) we see that there is a homotopy
\[h_\lambda: \alpha \simeq \beta~~\mathrm{in}~~\{u: \mathbb{I} \rightarrow \mathbb{S}^n: u(0)=u(1)\},\]
\[h_\lambda : \mathbb{S}^{n-1} \rightarrow \{u: \mathbb{I} \rightarrow \mathbb{S}^n: u(0)=u(1)=s_\lambda\} =: \Omega(\mathbb{S}^n,s_\lambda)\]
defined similarly like $\alpha$ and $\beta$ above. Of course,
there is a rotation $R$ of $\mathbb{S}^n$ such that
\[\alpha \simeq \bar{R}\circ\beta~~\mathrm{in}~~\Omega \mathbb{S}^n\]
with the homeomorphism $\bar{R} :\Omega(\mathbb{S}^n,s_1)
\rightarrow \Omega \mathbb{S}^n$, $\bar{R}(u)(t):= R(u(t))$.

\begin{Lemma} For every $n \geq 3$
\[\alpha_\star: H_{n-1}\mathbb{S}^{n-1} \rightarrow H_{n-1}(\Omega
\mathbb{S}^n\setminus \Omega_{\mathrm{tf}})\] is an isomorphism.
\end{Lemma}
Proof. Denote by $JX$ the James reduced product of $X$. Since the
pair $(J\mathbb{S}^{n-1},\mathbb{S}^{n-1})$ is $(2n-3)$-connected
\cite{P2} p.52, inclusion induces an isomorphism
\[H_q \mathbb{S}^{n-1} \rightarrow H_q J\mathbb{S}^{n-1}\]
for $q+1 \leq 2n-3$. By the James-Puppe Theorem $J\mathbb{S}^{n-1}
\simeq \Omega \mathbb{S}^n$ and the homomorphism
\[\beta_\star : H_{n-1}\mathbb{S}^{n-1} \rightarrow H_{n-1}\Omega(\mathbb{S}^n, s_1)\]
is an isomorphism \cite{H1}4.J, \cite{P1}. $\Box$\vspace{2mm}

(3.3) [Embeddings $\gamma_\omega$.] Fix any $\omega \in
\Omega_{\mathrm{tf}}$. For every $x \in \mathbb{S}^{n-1}$ denote
by $\mu(x)$ the meridian arc from $s_0$ to $x$, $\nu(x):=
\mu(-x)$. Then
\[\gamma_\omega(x):= \omega_{\mu(x)}\]
with its definition in 2. Note that embeddings $\alpha$ and
$\gamma_\omega$ are equivariant.

\section{Proof of Theorem 1.}
(A)~Fix a map $f: \Omega \mathbb{S}^n \setminus
\Omega_{\mathrm{tf}} \rightarrow \mathbb{R}^k$ with $k \leq n-1$.
Then
\[\{\omega_{\mu(x)}: \omega \in \Omega_{\mathrm{tf}}, x \in \mathbb{S}^{n-1}\} \approx \Omega_{\mathrm{tf}}\times \mathbb{S}^{n-1}.\]
If $P(Y,y_0)$ is the space of all paths in $Y$ from $y_0$ then
\[\Omega_{\mathrm{tf}} \approx P(\mathbb{S}^n,s_0)\supset P(\mathbb{S}^1\setminus \{-s_0\},s_0)\approx P(\mathbb{R},0)\supset \mathbb{R}^{d+1}
\supset \mathbb{S}^d,\] $P(\mathbb{R},0)$ is an
$\infty$-dimensional vector space. Briefly, $\mathbb{S}^d \subset
\Omega_{\mathrm{tf}}$ for every $d$. We apply Corollary 1 \cite{N}
to the following objects:
\[X := \{\omega_{\mu(x)}: \omega \in \mathbb{S}^d, x \in \mathbb{S}^{n-1}\} \approx \mathbb{S}^d\times \mathbb{S}^{n-1},\]
\[B = B^\prime := \mathbb{S}^d,~~ M:= \mathbb{S}^{n-1},~~ X^\prime:= \mathbb{S}^d\times\mathbb{R}^k,\]
\[T:X \rightarrow X,~~T(\omega,x)=(\omega,-x)~~\mathrm{i.e.}~~T(\omega_{\mu(x)}) = (\omega_{\mu(x)})^\star\]
and to the fibre-preserving map
\[X \ni (\omega,x) \stackrel{F}{\mapsto} (\omega, f(\omega_{\mu(x)}))\in X^\prime.\]
Put $C_F:= \{(\omega,x)\in X: F(\omega,x) = F(T(\omega,x))\}$.
Then
\[C_F = \{\omega_{\mu(x)}\in X: f(\omega_{\mu(x)}) = f((\omega_{\mu(x)})^\star)\} \subset A_f\]
and by Corollary 1 \cite{N},
\[H^{d+n-1-k}(C_F/T; \mathbb{Z}_2) \neq 0,~~H^{d+n-1-k}(C_F; \mathbb{Z}_2) \neq 0,\]
\[\mathrm{dim} A_f \geq \mathrm{dim} C_F \geq d+n-1-k \geq d~~\mathrm{for~every}~~d.\]

(B)~Suppose on the contrary that there is an equivariant map
\[\varphi: \Omega \mathbb{S}^n \setminus \Omega_{\mathrm{tf}}
\rightarrow \mathbb{S}^{n-1}.\] By Borsuk's Theorem the mapping
$\varphi \circ \gamma_{s_0}$ is essential. But
\[\gamma_{s_0}(x) = (s_0)_{\mu(x)} \simeq (s_0)_{\mu(s_1)} = \mathrm{const}\]
by homotopy
\[\left((s_0)_{(1-\lambda)\mu(x)}\right)_{\lambda \mu(s_1)},\]
a contradiction. $\Box$

\begin{Remark} Minor changes in the above proof show that Theorem
1 holds with $\mathbb{S}^n$ replaced by any $n$-dimensional
topological manifold.
\end{Remark}

\section{Proof of Proposition 1.}
Consider the Gysin exact sequence of the fibering $\mathbb{S}^0
\rightarrow E \rightarrow B$ of the form

\[\longrightarrow H^{q-2}(E;\mathbb{Z}_2) \longrightarrow H^{q-2}(B;\mathbb{Z}_2) \stackrel{\cup w_1}{\longrightarrow} H^{q-1}(B;\mathbb{Z}_2)
\stackrel{p^\star}{\longrightarrow} H^{q-1}(E;\mathbb{Z}_2)
\longrightarrow
\]
By Lemma 1 and \cite{H2} Example 1.5,
\[H_j E = H_j \Omega \mathbb{S}^n = \left\{\begin{array}{ll}
\mathbb{Z}& \mathrm{for}~~ j\equiv 0 ~(\mathrm{mod}~ (n-1))\\
0 & \mathrm{otherwise},\end{array}\right.\]
\[H^j(E;\mathbb{Z}_2) = \left\{\begin{array}{ll}
\mathbb{Z}_2& \mathrm{for}~~ j\equiv 0 ~(\mathrm{mod}~ (n-1))\\
0 & \mathrm{otherwise}.\end{array}\right.\] Since
$\pi_1E=\pi_2\mathbb{S}^n=0$,
$\pi_1B=\mathbb{Z}_2=H^1(B,\mathbb{Z}_2)$. By Gysin sequence,
\[H^j(B,\mathbb{Z}_2) = \mathbb{Z}_2~~\mathrm{for}~j \leq n-2~~\mathrm{and}~(w_1)^{n-1}\neq 0.\]

By Lemma 2 and the Universal Coefficient Theorem for cohomology
the homomorphism $\alpha^\star$ in the commutative diagram

\[\begin{array}{ccc}
H^{n-1}(E;\mathbb{Z}_2)& \stackrel{\alpha^\star}{\longrightarrow}& H^{n-1}(\mathbb{S}^{n-1},\mathbb{Z}_2)\\
p^\star\uparrow& & \uparrow \gamma^\star\\
H^{n-1}(B,\mathbb{Z}_2)&
\stackrel{\bar{\alpha}~^\star}{\longrightarrow}&
H^{n-1}(\mathbb{RP}^{n-1},\mathbb{Z}_2)
\end{array}
\]

is an isomorphism. Thus $p^\star =0$, because $\gamma^\star=0$.
Then $H^{n-1}(B;\mathbb{Z}_2)=\mathbb{Z}_2$ and by Gysin sequence,
$H^n(B;\mathbb{Z}_2)=0$.~ $\Box$

\begin{Remark} Further application of Gysin sequence yields
\[H^j(B;\mathbb{Z}_2)=0~~\mathrm{for}~~n <j<2n-2,\]
\[H^{2n-2}(B;\mathbb{Z}_2) = \mathbb{Z}_2~~\mathrm{with~a~generator}~~u,\]
\[H^j(B;\mathbb{Z}_2) = \mathbb{Z}_2~~\mathrm{for}~~2n-2 <j<3n-3,\]
\[H^{3n-3}(B;\mathbb{Z}_2) \ni u\cup(w_1)^{n-1} \neq 0.\]
It would be interesting to determine $H^j(B;\mathbb{Z}_2)$ for all
$j \geq 3n-3$.
\end{Remark}

\section{An example.} For $(f_1,\ldots, f_k) = f: \Omega \mathbb{R}^n \rightarrow
\mathbb{R}^k$ consider the equation
\[(\star)~~f(\alpha) = f(\alpha^\star),~~\mathrm{where}\]
\[f_j(\alpha) = \int_0^1 ||\alpha(t)-\beta_j(t)||^2
dt\] with fixed paths $\beta_j : \mathbb{I} \rightarrow
\mathbb{R}^n$, $j=1,\ldots,k$. The equation $(\star)$ is
equivalent to the system

\[(\star\star)~~\int_0^1 \langle\alpha(t), \beta_j(1-t) - \beta_j(t)\rangle dt =0,~~j=1,\ldots,k.\]

To see that $\mathrm{dim} A_f = \infty$ it suffices to specify
\[\alpha_x(t):= \left\{\begin{array}{lll}
x(t)& \mathrm{for}& t\in [0, \frac{1}{4}]\\
&&\\
x(\frac{1}{2}-t)& \mathrm{for}& t \in [\frac{1}{4},
\frac{1}{2}]\\
&&\\
0& \mathrm{for}& t \in [\frac{1}{2},1] \end{array}\right.\] for
any $x: [0,\frac{1}{4}] \rightarrow \mathbb{R}^n$, $x(0) = 0$, and
rewrite ($\star\star$) for $\alpha = \alpha_x$ in the form of a
system of $k$ homogeneous linear equations

\[\int_0^{\frac{1}{4}}\langle x(t), \beta_j(1-t) - \beta_j(t) + \beta_j(\frac{1}{2}+t) - \beta_j(\frac{1}{2}-t)\rangle dt = 0,~~j=1,\ldots,k; \]
in the $\infty$-dimensional space $C_0\left([0,\frac{1}{4}];
\mathbb{R}^n\right)$. Of course, $x \neq 0$ iff $\alpha_x \neq
(\alpha_x)^\star$.\vspace{2mm}

In this example $k$ and $n$ were any positive integers. The
natural question arises: does Theorem 1 hold for $k > n$?


\newpage
\begin{thebibliography}{999}
\bibitem{B} K. Borsuk, {\it Drei S\"{a}tze \"{u}ber die $n$-dimensionale euklidische Sph\"{a}re}, Fund. Math. 20 (1933), 177--190.
\bibitem{H1} A. Hatcher, Algebraic Topology, Cambridge University
Press, 2002.
\bibitem{H2} A. Hatcher, www.math.cornell.edu/~hatcher/SSAT/SSch1.pdf
\bibitem{J} J. Jaworowski, {\it A continuous version of
Borsuk-Ulam theorem}, Proc. Amer. Math. Soc. 82 (1981), 112--114.
\bibitem{N} M. Nakaoka, {\it Equivariant point theorems for
fibre-preserving maps},\\ Osaka J. Math. 21(1984), 809--815.
\bibitem{P1} M. M. Postnikov, {\it Lekcii po algebrai\v{c}eskoj
topologii: osnovy teorii gomotopij}~\footnote{Lectures in
algebraic topology: elements of homotopy theory.}, Moskva, 1984.
\bibitem{P2} M. M. Postnikov, {\it Lekcii po algebrai\v{c}eskoj
topologii: teori\^{a} gomotopij kleto\v{c}nyh
prostranstv}~\footnote{Lectures in algebraic topology: homotopy
theory of cell complexes.}, Moskva, 1985.
\end{thebibliography}
\end{document}